\documentclass[12pt]{article}

\usepackage{amssymb}
\usepackage[mathscr]{euscript}

\def\<{\langle} \def\>{\rangle}

\title{On Konopelchenko's representation formula for surfaces in 4 dimensions}
\author{Fr\'ed\'eric H\'elein}

\date{April 9, 2001}

\begin{document}
\maketitle
In \cite{HR2} Pascal Romon and I proposed a Weierstrass type formula
for conformal parametrisations of Lagrangian surfaces in $\Bbb{R}^4$,
which relies on a kind of Dirac equation. Similarities with the
Weierstrass formula found by B.G. Konopelchenko for surfaces in $\Bbb{R}^3$
\cite{K1} were observed. Recently Franz Pedit pointed out to me that
B.G. Konopelchenko did generalize his formula for representing surfaces in
$\Bbb{R}^4$ \cite{K2}. The purpose of this note is simply to stress the fact
that our representation formula for Lagrangian surfaces can be recovered from
Konopelchenko's one. As a byproduct we can deduce a simplification of
Konopelchenko's representation by using quaternions.

\section{Weierstrass representations of Lagrangian surfaces}
Here we recall results in \cite{HR2}. We identify $\Bbb{R}^4$ with $\Bbb{C}^2$,
with the Hermitian scalar product $\langle v,w,\rangle _H:=
v^1\overline{w^1} + v^2\overline{w^2}$. Hence, as a four-dimensional space it
has the Euclidean scalar product $\langle .,.\rangle _E$ and the symplectic
form $\omega$, such that $\langle .,.\rangle _H = \langle .,.\rangle _E
-i\omega$. Let $\Omega$ be an open simply connected domain of $\Bbb{R}^2$.
Then all smooth conformal Lagrangian immersions of $\Omega$ can be obtained as
follows. Start with any smooth function $\beta:\Omega\longrightarrow \Bbb{R}$
and denote

$$p := {1\over 2}{\partial \beta\over \partial \overline{z}}.$$
Then take any pair of functions $s_1,s_2:\Omega \longrightarrow \Bbb{C}$ which
are solutions of

\begin{equation}\label{1.1}
\left\{
\begin{array}{ccl}
\displaystyle 
{\partial s_1\over \partial \overline{z}} & = & -\overline{p}\overline{s_2}\\
&&\\
\displaystyle {\partial \overline{s_2}\over \partial z} & = & ps_1.
\end{array}\right.
\end{equation}
Then a smooth conformal Lagrangian immersion $Y:\Omega\longrightarrow
\Bbb{C}^2$ is given by

\begin{equation}\label{1.2}
\begin{array}{ccl}
Y^1(z)+iY^2(z) & = & \displaystyle C^1+iC^2 + \int_{\Gamma}{1\over
2}e^{i{\beta\over 2}}\left[ (s_1+is_2)dz +
(\overline{s_1}+i\overline{s_2})d\overline{z}\right]\\
Y^3(z)+iY^4(z) & = &
\displaystyle C^3+iC^4 + \int_{\Gamma}{1\over 2}e^{i{\beta\over 2}}\left[
(-is_1-s_2)dz + (i\overline{s_1}+\overline{s_2})d\overline{z}\right],
\end{array} \end{equation}
where $C^i$ are constants of integration and $\Gamma$ is any based path with
endpoint $z$. The existence of these formulas hence means that the 1-forms
in the integrals are closed ones. The function $\beta$ has an important
geometrical interpretation in symplectic geometry, it is called the
Lagrangian angle (see \cite{HR1}, \cite{SW}).\\

In the following we shall use these formulas with other variables: 
set $(X^1,X^2,X^3,X^4):= (Y^1,Y^3,Y^2,Y^4)$ and drop the constants $C^i$.
Then (\ref{1.2}) reads

\begin{equation}\label{1.3}
\begin{array}{ccl}
X^1(z)+iX^2(z) & = &  \displaystyle \int_{\Gamma}
s_1\cos {\beta\over 2}dz -\overline{s_2}\sin {\beta\over 2}d\overline{z}\\
X^3(z)+iX^4(z) & = &  \displaystyle \int_{\Gamma}
s_1\sin {\beta\over 2}dz + \overline{s_2}\cos {\beta\over 2}d\overline{z}.
\end{array}
\end{equation}

\section{Representation of surfaces in $\Bbb{R}^4$}
The following is due to Konopelchenko \cite{K2} (see also \cite{p1} and \cite{p2}
for more results). We change slightly the
notations (for comparaison, Konopelchenko uses variables $\psi_1$, $\psi_2$,
$\phi_1$ and $\phi_2$; here we have $s_1=\phi_1$, $s_2=\overline{\psi_1}$,
$t_1=\phi_2$ and $t_2=\overline{\psi_2}$. Moreover we have changed the sign
of $X^1$ and $X^2$). Let $p:\Omega\longrightarrow \Bbb{C}$ be any smooth
function. Assume that $s_1, s_2,t_1,t_2:\longrightarrow \Bbb{C}$ are solutions
of

\begin{equation}\label{2.1}
\left\{
\begin{array}{ccl}
\displaystyle {\partial s_1\over \partial \overline{z}} & = & -\overline{p}\overline{s_2}\\
&&\\
\displaystyle {\partial \overline{s_2}\over \partial z} & = & ps_1
\end{array}\right.
\end{equation}
and
\begin{equation}\label{2.2}
\left\{
\begin{array}{ccl}
\displaystyle {\partial t_1\over \partial \overline{z}} & = & -p\overline{t_2}\\
&&\\
\displaystyle {\partial \overline{t_2}\over \partial z} & = & \overline{p}t_1.
\end{array}\right.
\end{equation}
Then the 1-forms $s_1t_1dz - \overline{s_2}\overline{t_2}d\overline{z}$ and
$s_1t_2dz + \overline{s_2}\overline{t_1}d\overline{z}$ are
closed and the functions given by

\begin{equation}\label{2.3}
\begin{array}{ccl}
X^1(z)+iX^2(z) & = &  \displaystyle \int_{\Gamma}
s_1t_1dz - \overline{s_2}\overline{t_2}d\overline{z}\\
X^3(z)+iX^4(z) & = &  \displaystyle \int_{\Gamma}
s_1t_2dz + \overline{s_2}\overline{t_1}d\overline{z}
\end{array}
\end{equation}
define a conformal immersion into $\Bbb{R}^4$.

Here we see that if we set $t_1= \cos {\beta\over 2}$ and $t_2=\sin {\beta\over
2}$, then (\ref{2.3}) implies (\ref{1.3}). We can check that actually these
datas are a simple solution to (\ref{2.2}) for $p = {1\over 2}{\partial
\beta\over \partial \overline{z}}$.

\section{A formulation using quaternions}
Inspired by \cite{HR2}, where an alternative presentation using quaternions was
given, we can reformulate (\ref{2.1}), (\ref{2.2}) and (\ref{2.3}). We denote
$\Bbb{H}$ the set of quaternions and $i,j,k$ its three imaginary number. We
now view any function $f:\Omega\longrightarrow \Bbb{R}^4$ as a function into
$\Bbb{H}$ by denoting

$$f = f^1 + if^2 + jf^3 + kf^4 = (f^1 + if^2) + (f^3+if^4)j.$$
And we think $dz := dx +idy$ and $\overline{dz}:= dx-idy$, as 1-forms with
values in $\Bbb{H}$. We can hence define left and right derivatives by

$$df = (\partial f/\partial z)dz + (\partial f/\partial
\overline{z})d\overline{z} = dz(\partial z\setminus \partial f) +
d\overline{z}(\partial \overline{z}\setminus \partial f).$$

Setting $a:= s_1+js_2$ and $b:= t_1+t_2j$, we compute that
$(\partial a/\partial \overline{z}) =  {\partial s_1\over \partial
\overline{z}} + j {\partial s_2\over \partial \overline{z}}$ and
$(\partial \overline{z}\setminus \partial b) = {\partial t_1\over \partial
\overline{z}} + {\partial t_2\over \partial \overline{z}}j$. Hence
(\ref{2.1}) and (\ref{2.2}) are translated into

\begin{equation}\label{3.1}
(\partial a/\partial \overline{z}) = apj\quad \hbox{and} \quad 
(\partial \overline{z}\setminus \partial b) = pjb.
\end{equation}
Now (\ref{2.3}) reads

\begin{equation}\label{3.2}
X^1+iX^2+jX^3+kX^4 = \int_{\Gamma}
a\,dz\,b.
\end{equation}
And it is straightforward to prove that (\ref{3.1}) is the
integrability condition for (\ref{3.2}), since \footnote{beware that
$d(dz\;f) = -dz\wedge d\overline{z}(\partial \overline{z}\setminus \partial f)$}

$$\begin{array}{ccl}
d(a\,dz\,b) & = & (\partial a/\partial \overline{z})d\overline{z}\wedge dz\,b
- a\,dz\wedge d\overline{z}(\partial \overline{z}\setminus \partial b)\\
 & = & 2(apj\,i\,b + a\,i\,pjb) dx\wedge dy\\
 & = & 2ap(ji + ij)b\,dx\wedge dy = 0.
\end{array}$$

\noindent {\bf Aknowledgements} I am pleased to thank Franz Pedit for valuable
discussions.

\noindent {\em Fr\'ed\'eric H\'elein,
CMLA, ENS de Cachan, 61, avenue du Pr\'esident Wilson, 94235 Cachan Cedex,
France}, helein@cmla.ens-cachan.fr

\end{document}